# Verdier Stratifications and ($w_f$)-Stratifications in o-minimal Structures

## Ta Lê Loi

**Abstract** – We prove the existence of Verdier stratifications for sets definable in any o-minimal structure on $(\mathbf{R}, +, \cdot)$. It is also shown that the Verdier condition (w) implies the Whitney condition (b) in o-minimal structures on $(\mathbf{R}, +, \cdot)$. As a consequence the Whitney Stratification Theorem holds. The existence of ($w_f$)-stratification of functions definable in polynomially bounded o-minimal structures is presented.

**0. Introduction.**

**0.1 Definition.** An *o-minimal structure* on the real field $(\mathbf{R}, +, \cdot)$ is a family $\mathcal{D} = (\mathcal{D}_n)_{n \in \mathbf{N}}$ such that for each $n \in \mathbf{N}$:

(1) $\mathcal{D}_n$ is a boolean algebra of subsets of $\mathbf{R}^n$.
(2) If $A \in \mathcal{D}_n$, then $A \times \mathbf{R}$ and $\mathbf{R} \times A \in \mathcal{D}_{n+m}$.
(3) If $A \in \mathcal{D}_{n+1}$, then $\pi(A) \in \mathcal{D}_n$, where $\pi : \mathbf{R}^{n+1} \longrightarrow \mathbf{R}^n$ is the projection on the first $n$ coordinates.
(4) $\mathcal{D}_n$ contains $\{x \in \mathbf{R}^n : P(x) = 0\}$ for all polynomials $P \in \mathbf{R}[X_1, \cdots, X_n]$.
(5) Each set belonging to $\mathcal{D}_1$ is a finite union of intervals and points. (o- minimality)

A set belonging to $\mathcal{D}$ is called *definable* (in this structure). *Definable maps* are maps whose graphs belonging to $\mathcal{D}$.

Many results in Semialgebraic Geometry and Subanalytic Geometry hold true in the theory of o-minimal structures on the real field. Recently, o-minimality of many interesting structrures on $(\mathbf{R}, +, \cdot)$ has been established, for example, structures generated by the exponential function [W1](see also [LR] and [DM1]), real power functions [M2], Pfaffian functions [W2] or restricted Gevrey functions [DS]. For more details on o-minimal structures we refer the readers to [D] and [DM] (compare with [S]).

**0.2** We now outline the main results of this paper. Let $\mathcal{D}$ be an o-minimal structure on $(\mathbf{R}, +, \cdot)$. In section 1, we prove that $\mathcal{D}$ admits Verdier Stratification. We also show that the Verdier condition (w) implies the Whitney condition (b) in $\mathcal{D}$. Thus, Whitney Stratification Theorem holds true in $\mathcal{D}$. These improve results in [L1] (see also [DM]). Note that the theorems were proved for subanalytic sets in [V] and [ŁSW] (see also [DW]), the former based on Hironaka's Desingularization, and the latter on

---

[0] *Key words:* o-minimal structures - Verdier Stratification - ($w_f$)-Stratification
[0] *Mathematics Subject Classification:* 32S60 14P10 14B5.

Puiseux's Theorem. But, in general, neither tools can be applied to sets belonging to o-minimal structures (e.g. to the set $\{(x,y) \in \mathbf{R}^2 : y = \exp(-1/x), x > 0\}$ in the structure generated by the exponential function). Section 2 is devoted to the study of stratifications of definable functions. In general, definable functions cannot be stratified to satisfy the strict Thom condition ($w_f$). However, if $\mathcal{D}$ is polynomially bounded, then it admits ($w_f$)-stratification. Our proof of this assertion is based on piecewise uniform asymptotics for definable functions, instead of Pawłucki's version of Puiseux's theorem with parameters, that is used in [KP] to prove the assertion for subanalytic functions.

**0.3 Notation.** Throughout this paper, let $\mathcal{D}$ denote some fixed, but arbitrary, o-minimal structure on $(\mathbf{R}, +, .)$. "*Definable*" means definable in $\mathcal{D}$. If $\mathbf{R}^k \times \mathbf{R}^l \ni (y,t) \mapsto f(y,t) \in \mathbf{R}^m$ is a differentiable function, then $D_1 f$ denotes the derivative of $f$ with respect to the first variables $y$. As usual, $d(\cdot,\cdot), \|\cdot\|$ denote the Euclidean distance and norm respectively. Besides, Cell Decomposition [DM. Th. 4.2], and Definable Choice [DM. Th. 4.5] will be often referred in our arguments without the citations.

**Acknowledgements.** I wish to thank the Fields Institute, University of Toronto, where this paper was written, for hospitality and support. I also thank to Chris Miller for many helpful suggestions.

**1. Verdier Stratifications**

**1.1 Verdier condition.** Let $\Gamma, \Gamma'$ be $C^1$ submanifolds of $\mathbf{R}^n$ such that $\Gamma \subset \overline{\Gamma}' \setminus \Gamma'$. Let $y_0$ be a point of $\Gamma$. We say that the pair $(\Gamma, \Gamma')$ satisfies the Verdier condition at $y_0$ if the following holds:

(w) There exists a constant $C > 0$ and a neighborhood $U$ of $y_0$ in $\mathbf{R}^n$ such that

$$\delta(T_y\Gamma, T_x\Gamma') \leq C\|x - y\| \quad \text{for all } x \in \Gamma' \cap U, y \in \Gamma \cap U,$$

where $T_y\Gamma$ denotes the tangent space of $\Gamma$ at $y$, and $\delta(T, T') = \sup_{v \in T, \|v\|=1} d(v, T')$ is the distance of vector subspaces of $\mathbf{R}^n$.
Note that (w) is invariant under $C^2$-diffeomorphisms.

**1.2 Definition.** Let $p$ be a positive integer. A $C^p$ *stratification* of $\mathbf{R}^n$ is a partition $\mathcal{S}$ of $\mathbf{R}^n$ into finitely many subsets, called strata, such that:
(S1) Each stratum is a connected $C^p$ submanifold of $\mathbf{R}^n$ and also a definable set.
(S2) For every $\Gamma \in \mathcal{S}$, $\overline{\Gamma} \setminus \Gamma$ is a union of some of the strata.

We say that $\mathcal{S}$ is *compatible with* a class $\mathcal{A}$ of subsets of $\mathbf{R}^n$ if each $A \in \mathcal{A}$ is a finite union of some strata in $\mathcal{S}$.
A $C^p$ *Verdier stratification* is a $C^p$ stratification $\mathcal{S}$ such that for all $\Gamma, \Gamma' \in \mathcal{S}$, if $\Gamma \subset \overline{\Gamma'} \setminus \Gamma'$, then $(\Gamma, \Gamma')$ satisfies the condition (w) at each point of $\Gamma$.

**1.3 Theorem** (Verdier Stratification). *Let $p$ be a positive integer. Then given definable sets $A_1, \cdots, A_k$ contained in $\mathbf{R}^n$, there exists a $C^p$ Verdier stratification of $\mathbf{R}^n$ compatible with $\{A_1, \cdots, A_k\}$.*

We first make an observation similar to that of [ŁSW. Prop. 2].
Let (P) be a (local) property for the pairs $(\Gamma, \Gamma')$ at $y$ in $\Gamma$, where $\Gamma, \Gamma'$ being subsets of $\mathbf{R}^n$. Put $P(\Gamma, \Gamma') = \{y \in \Gamma : (\Gamma, \Gamma') \text{ satisfies (P) at } y\}$.

**1.4 Proposition.** *Suppose that for every pair $(\Gamma, \Gamma')$ of definable $C^p$ submanifolds of $\mathbf{R}^n$ with $\Gamma \subset \overline{\Gamma'} \setminus \Gamma'$, the set $P(\Gamma, \Gamma')$ is definable and $\dim(\Gamma \setminus P(\Gamma, \Gamma')) < \dim \Gamma$. Then given definable sets $A_1, \cdots, A_k$ contained in $\mathbf{R}^n$, there exists a $C^p$ stratification $\mathcal{S}$ of $\mathbf{R}^n$ compatible with $\{A_1, \cdots, A_k\}$ and has the following*

(P)  $P(\Gamma, \Gamma') = \Gamma$  *for all* $\Gamma, \Gamma' \in \mathcal{S}$ *with* $\Gamma \subset \overline{\Gamma'} \setminus \Gamma'$.

*Proof.* We can construct, by decreasing induction on $d \in \{0, \cdots, n\}$, partitions $\mathcal{S}^d$ of $\mathbf{R}^n$ into $C^p$-cells compatible with $\{A_1, \cdots, A_k\}$, such that $\mathcal{S}^d$ has properties (S1)(S2) and the following property:

$(P_d)$ $\quad P(\Gamma, \Gamma') = \Gamma$  for all $\Gamma, \Gamma' \in \mathcal{S}^d$ with $\Gamma \subset \overline{\Gamma'} \setminus \Gamma'$ and $\dim \Gamma \geq d$.

Indeed, by Cell Decomposition and the fact that $\dim(\overline{A} \setminus A) < \dim A$, for all definable set $A$, we can construct a $C^p$ cell decomposition of $\mathbf{R}^n$ compatible with $\{A_1, \cdots, A_k\}$ and has (S1)(S2). This cell decomposition can be refined to satisfy $(P_d)$ by the assumption.
Obviously, $\mathcal{S} = \mathcal{S}^0$ is a desired stratification. $\quad \diamond$

By the proposition, Theorem 1.3 follows from the the following.

**1.5 Theorem.** *Let $\Gamma, \Gamma'$ be definable sets and $C^p$-submanifolds of $\mathbf{R}^n$. Suppose that $\Gamma \subset \overline{\Gamma'} \setminus \Gamma'$. Then $W = \{y \in \Gamma : (\Gamma, \Gamma') \text{ satisfies (w) at } y\}$ is definable, and $\dim(\Gamma \setminus W) < \dim \Gamma$.*

To prove Theorem 1.5 we prepare some lemmas.

**1.6 Lemma.** *Under the notation of Theorem 1.5, $W$ is a definable set.*

*Proof.* Note that the Grassmanian $G_k(\mathbf{R}^n)$ of $k$-dimensional linear subspaces of $\mathbf{R}^n$ is semialgebraic, and hence definable; $\delta$ and the tangent map: $\Gamma \ni x \mapsto T_x\Gamma \in G_{\dim \Gamma}(\mathbf{R}^n)$ are also definable. Therefore,

$$W = \{y_0 : \quad y_0 \in \Gamma, \exists C > 0, \exists t > 0, \forall x \in \Gamma', \forall y \in \Gamma \\ (\|x - y_0\| < t, \|y - y_0\| < t \Rightarrow \delta(T_y\Gamma, T_x\Gamma') \leq C\|x - y\|)\}$$

is a definable set. $\quad \diamond$

**1.7 Lemma** (Wing Lemma). *Let $V \subset \mathbf{R}^k$ be open and definable, and $S \subset \mathbf{R}^k \times \mathbf{R}^l$ be definable. Suppose $V \subset \overline{S} \setminus S$. Then there exists an open subset $U$ of $V$, $\alpha > 0$, and a definable map $\bar{\rho} : U \times (0, \alpha) \longrightarrow S$, of class $C^p$, such that $\bar{\rho}(y, t) = (y, \rho(y, t))$ and $\|\rho(y, t)\| = t$, for all $y \in U, t \in (0, \alpha)$.*

*Proof.*(c.f. [L1. Lemma 2.7]) Let

$$A = \{(y, x, t) : y \in V, x \in S, 0 < t < 1, \|x - y\| < t, \pi(x) = y\},$$

where $\pi : \mathbf{R}^k \times \mathbf{R}^l \longrightarrow \mathbf{R}^k$ is the natural projection.

Note that $A$ is definable. Let $\pi_2(y, x, t) = (y, t)$, and $\pi_2(A)_y = \{t : (t, y) \in \pi_2(A)\}$. Define

$$\epsilon(y) = \inf \pi_2(A)_y , \quad (\inf \emptyset := 1).$$

Then $\epsilon : V \longrightarrow \mathbf{R}$ is definable, and if $\epsilon(y) > 0$ then $(0, \epsilon(y)) \cap \pi_2(A)_y = \emptyset$.

*Claim:* $\dim\{y \in V : \epsilon(y) > 0\} < \dim V = k$.

Suppose to the contrary that the dimension is $k$. Then, by Cell Decomposition, there is an open ball $B \subset V$ and $c > 0$ such that $\epsilon > c$ on $B$. This implies $B \not\subset \overline{V} \setminus S$, a contradiction.

Now let $V_0 = \{y \in V : \epsilon(y) = 0\}$. Then $\dim V_0 = k$, and, by the definition, $V_0 \times (0, 1) \subset \pi_2(A)$. By Definable Choice and Cell Decomposition, there exists an open set $V' \subset V_0$, $\delta > 0$, and a continuous definable map: $V' \times (0, \delta) \longrightarrow S$, $(y, t) \mapsto (y, \theta(y, t))$. Let $\tau(y) = \sup_{0<s<\delta} \|\theta(y, s)\|$. Then for $y \in V'$, $t < \tau(y)$, there exists $x \in S$, such that $\pi(x) = y$ and $\|x - y\| = t$. Again by Definable Choice and Cell Decomposition it is easy to prove the existence of the $U, \alpha, \bar{\rho}$ satisfying the demands of the lemma. ◇

To control the tangent spaces we need the following lemma.

**1.8 Lemma.** *Let $U \subset \mathbf{R}^k$ be an open definable set, and $M : U \times (0, \alpha) \longrightarrow \mathbf{R}^l$ be a $C^1$ definable map. Suppose there exists $K > 0$ such that $\|M(y, t)\| \leq K$, for all $y \in U$ and $t \in (0, \alpha)$. Then there exists a definable set $F$, closed in $U$ with $\dim F < \dim U$, and continuous definable functions $C, \tau : U \setminus F \longrightarrow \mathbf{R}_+$, such that for all $y$ in $U \setminus F$*

$$\|D_1 M(y, t)\| \leq C(y) , \quad \text{for all } t \in (0, \tau(y)).$$

*Proof.* It suffices to prove for $l = 1$. Suppose the assertion of the lemma is false. Since $\{y \in U : \lim_{t \to 0^+} \|D_1 M(y, t)\| = +\infty\}$ is definable, there is an open subset $B$ of $U$, such that

$$\lim_{t \to 0^+} \|D_1 M(y, t)\| = +\infty, \quad \text{for all } y \text{ in } B.$$

By monotonicity [DM. Th. 4.1], for each $y \in B$, there is $s > 0$ such that $t \mapsto \|D_1 M(y, t)\|$ is strictly decreasing on $(0, s)$. Let

$$\tau(y) = \sup\{s : \|D_1 M(y, \cdot)\| \text{ is strictly decreasing on } (0, s)\}.$$

Note that $\tau$ is a definable function, and, by Cell Decomposition, $\tau$ is continuous on an open subset $B'$ of $B$, and $\tau > \alpha'$ on $B'$, for some $\alpha' > 0$. Let $\psi(t) = \inf\{\|D_1 M(y, t)\| : y \in B', 0 < t < \alpha'\}$. Shrinking $B'$, we can assume that $\lim_{t \to 0^+} \psi(t) = +\infty$. Then, for each $y \in B'$, we have

$$\|D_1 M(y, t)\| > \psi(t) , \quad \text{for all } t \in (0, \alpha').$$

This implies $|M(y,t) - M(y',t)| > \psi(t)\|y - y'\|$, for all $y, y' \in B'$, and $t < \alpha'$.

Therefore, $\psi(t) \leq \dfrac{2K}{diam B'}$, for all $t \in (0, \alpha')$, a contradiction. ◇

**1.9 Proof of Theorem 1.5.** The first part of the theorem was proved in Lemma 1.6. To prove the second part we suppose, contrary to the assertion, that $\dim(\Gamma \setminus W) = \dim \Gamma = k$.

Since (w) is a local property and invariant under $C^2$ local diffeomorphisms, we can suppose $\Gamma$ is an open subset of $\mathbf{R}^k \subset \mathbf{R}^k \times \mathbf{R}^{n-k}$. In this case $T_y\Gamma = \mathbf{R}^k$, for all $y \in \Gamma$. Then by the assumption, applying Lemma 1.7, we get an open subset $U$ of $\Gamma$, a $C^p$ definable map $\bar\rho : U \times (0, \alpha) \longrightarrow \Gamma'$ such that $\bar\rho(y,t) = (y, \rho(y,t))$ and $\|\rho(y,t)\| = t$, and, moreover, for each $y \in U$

$$\frac{\delta(\mathbf{R}^k, T_{(y,\rho(y,t))}\Gamma')}{\|\rho(y,t)\|} \to +\infty, \quad \text{when } t \to 0^+.$$

On the other hand, applying Lemma 1.8 to $M(y,t) := \dfrac{\rho(y,t)}{t}$, reducing $U$ and $\alpha$, we have

$$\|D_1\rho(y,t)\| \leq Ct, \quad \text{for all } y \in U, t \in (0, \alpha),$$

with some $C > 0$.

Note that $T_{(y,\rho(y,t))}\Gamma' \supset \operatorname{graph} D_1\rho(y,t)$. Therefore,

$$\frac{\delta(\mathbf{R}^k, T_{(y,\rho(y,t))}\Gamma')}{\|\rho(y,t)\|} \leq \frac{\|D_1\rho(y,t)\|}{\|\rho(y,t)\|} \leq C, \quad \text{for } y \in U, 0 < t < \alpha.$$

This is a contradiction. ◇

Note that Whitney condition (b) (defined in [Wh]) does not imply condition (w), even for algebraic sets (see [BT]). And, in general, we do not have (w) $\Rightarrow$ (b) (e.g. $\Gamma = (0,0)$, $\Gamma' = \{(x,y) \in \mathbf{R}^2 : x = r\cos r, y = r\sin r, r > 0\}$, or $\Gamma' = \{(x,y) \in \mathbf{R}^2 : y = x\sin(1/x), x > 0\}$). In o-minimal structures, by the finiteness, such spiral phenomena or oscillation cannot occur. The following is a version of Kuo-Verdier's Theorem (see [K] [V]).

**1.10 Proposition.** *Let $\Gamma, \Gamma' \subset \mathbf{R}^n$ be definable $C^p$-submanifolds ($p \geq 2$), with $\Gamma \subset \overline{\Gamma'} \setminus \Gamma'$. If $(\Gamma, \Gamma')$ satisfies the condition (w) at $y \in \Gamma$, then it satisfies the Whitney condition (b) at $y$.*

*Proof.* Our proof is an adaptation of [V. Theorem 1.5] and based on the following observation:

If $f : (0, \alpha) \longrightarrow \mathbf{R}$ is definable with $f(t) \neq 0$, for all $t$, and $\lim\limits_{t \to 0^+} f(t) = 0$, then, by Cell Decomposition and monotonicity [DM. Th.4.1], there is $0 < \alpha' < \alpha$, such that $f$ is of class $C^1$ and strictly monotone on $(0, \alpha')$. By Mean Value Theorem and Definable Choice, there exists a definable function $\theta : (0, \alpha') \to (0, \alpha')$ with $0 < \theta(t) < t$, such that $f(t) = f'(\theta(t))t$. Since $|f(t)| > |f(\theta(t))|$, by monotonicity, $\lim\limits_{t \to 0^+} \dfrac{f(t)}{f'(t)} = 0$.

Now we prove the theorem. By a $C^2$ change of local coordinates, we can suppose $\Gamma$ is an open subset of $\mathbf{R}^k \subset \mathbf{R}^k \times \mathbf{R}^l$ ($l = n - k$), and $y = 0$. Let $\pi : \mathbf{R}^k \times \mathbf{R}^l \longrightarrow \mathbf{R}^l$ be the orthogonal projection. Since $(\Gamma, \Gamma')$ satisfies (w) at 0, there exists $C > 0$ and a neigborhood $U$ of 0 in $\mathbf{R}^n$, such that

$$(*) \qquad \delta(T_y\Gamma, T_x\Gamma') \leq C\|x - y\|, \quad \text{for all } x \in \Gamma' \cap U, y \in \Gamma \cap U.$$

If the condition (b) is not satisfied at 0 for $(\Gamma, \Gamma')$, then there exists $\epsilon > 0$, such that $0 \in \overline{S} \setminus S$, where
$$S = \{x \in \Gamma' : \delta(\mathbf{R}\pi(x), T_x\Gamma') \geq 2\epsilon\}.$$

Since $S \cap \{x : \|x\| \leq t\} \neq \emptyset$, for all $t > 0$, by Curve selection [DM. Th.4.6], there exists a definable curve $\varphi : (0, \alpha) \longrightarrow S$, such that $\|\varphi(t)\| \leq t$, for all $t$. By the above observation, we can assume $\varphi$ is of class $C^1$. Write $\varphi(t) = (a(t), b(t)) \in \mathbf{R}^k \times \mathbf{R}^l$. Then $\|b'(t)\|$ is bounded. Since $\varphi((0, \alpha)) \subset \Gamma'$, $a \not\equiv 0$. Reducing $\alpha$, we can assume $a'(t) \neq 0$, for all $t$. Since $\lim_{t \to 0^+} a'(t)$ exists, we have $\delta(\mathbf{R}a'(t), \mathbf{R}a(t)) \to 0$, when $t \to 0$. Therefore

$$(**) \qquad \delta(\mathbf{R}a'(t), T_{\varphi(t)}\Gamma') \geq \epsilon, \quad \text{for all } t \text{ sufficiently small}.$$

On the other hand, we have

$$\delta(\mathbf{R}a'(t), T_{\varphi(t)}\Gamma') = \frac{1}{\|a'(t)\|}\delta(a'(t), T_{\varphi(t)}\Gamma') = \frac{1}{\|a'(t)\|}\delta(b'(t), T_{\varphi(t)}\Gamma')$$
$$\leq \frac{\|b'(t)\|}{\|a'(t)\|}\delta(\mathbf{R}b'(t), T_{\varphi(t)}\Gamma').$$

From $(*)$ and $(**)$, we have $\epsilon \leq C\|a(t)\|\dfrac{\|b'(t)\|}{\|a'(t)\|}$.

By the observation, the right-hand side of the inequality tends to 0 (when $t \to 0$), which is a contradiction. $\diamond$

From Theorem 1.3 and Proposition 1.10 we have

**1.11 Corollary.** *Whitney Stratification Theorem holds true in any o-minimal structure on the real field.*

## 2. ($w_f$)-Stratifications

Thoughout this section, let $X \subset \mathbf{R}^n$ be a definable set and $f : X \longrightarrow \mathbf{R}$ be a continuous definable function. Let $p$ be a positive integer.

**2.1 Definition.** A $C^p$ *stratification of* $f$ is a $C^p$ stratification $\mathcal{S}$ of $X$, such that for every stratum $\Gamma \in \mathcal{S}$, the restriction $f|_\Gamma$ is $C^p$ and of constant rank.
For each $x \in \Gamma$, $T_{x,f}$ denotes the tangent space of the level of $f|_\Gamma$ at $x$, i.e. $T_{x,f} = \ker d(f|_\Gamma)(x)$.
Let $\Gamma, \Gamma' \in \mathcal{S}$ with $\Gamma \subset \overline{\Gamma'} \setminus \Gamma'$. We say that the pair $(\Gamma, \Gamma')$ satisfies the *Thom condition* $(a_f)$ at $y_0 \in \Gamma$ if and only if the following holds:

($a_f$) for every sequence $(x_k)$ in $\Gamma'$, converging to $y_0$, we have

$$\delta(T_{y_0,f}, T_{x_k,f}) \longrightarrow 0 .$$

We say that $(\Gamma, \Gamma')$ satisfies the *strict Thom condition* $(w_f)$ at $y_0$ if:

$(w_f)$ there exists a constant $C > 0$ and a neightborhood $U$ of $y_0$ in $\mathbf{R}^n$, such that

$$\delta(T_{y,f}, T_{x,f}) \leq C\|x - y\| \quad \text{for all } x \in \Gamma' \cap U, y \in \Gamma \cap U.$$

Note that the conditions are $C^2$-invariant.

The existence of stratifications satisfying $(w_f)$ (and hence $(a_f)$) for subanalytic function was proved in [KP] (see also [B][KR]). For functions definable in o-minimal structures on the real field we have:

**2.2 Theorem.** *There exists a $C^p$ stratification of $f$ satisfying the Thom condition $(a_f)$ at every point of the strata.*

*Proof:* see [L2]  ◇

**2.3 Remark.** In general, definable functions cannot be stratified to satisfy the condition $(w_f)$. The following example is given by Kurdyka.
Let $f : [a,b] \times [0,+\infty) \longrightarrow \mathbf{R}$ be defined by $f(x,y) = y^x$ $(0 < a < b)$. Let $\Gamma = [a,b] \times 0$, and $\Gamma' = [a,b] \times (0,+\infty)$. Then the fiber of $f|_{\Gamma'}$ over $c \in \mathbf{R}_+$:

$$\{(x, y(x) = \exp(-\frac{1}{tx})) : x \in [a,b]\} , \quad t = -\frac{1}{\ln c}.$$

Then $\dfrac{y'(x)}{y(x)} = \dfrac{1}{tx^2} \to +\infty$, when $t \to 0^+$, for all $x \in [a,b]$,

i.e. $\dfrac{\delta(T_{x,f}, T_{(x,y(x)),f})}{\|y(x)\|}$ can not be locally bounded along $\Gamma$.

The remainder of this section is devoted to the proof of the existence of $(w_f)$-stratification of functions definable in polynomially bounded o-minimal structures.

**2.4 Definition.** A structure $\mathcal{D}$ on the real field $(\mathbf{R},+,.)$ is *polynomially bounded* if for every function $f : \mathbf{R} \longrightarrow \mathbf{R}$ definable in $\mathcal{D}$, there exists $N \in \mathbf{N}$, such that

$$|f(t)| \leq t^N , \quad \text{for all sufficiently large } t.$$

For example, the structure of global subanalytic sets, the structure generated by real power functions [M2], or by restricted Gevrey functions [DS] are polynomially bounded.

**2.5 Theorem.** *Suppose that $\mathcal{D}$ is polynomially bounded. Then there exists a $C^p$ stratification of $f$ satisfying the condition $(w_f)$ at each point of the strata.*

*Note* - The converse of the theorem is also true: If $\mathcal{D}$ is not polynomially bounded, then it must contain the exponential function, by [M1]. So the function given in Remark 2.3 is definable in $\mathcal{D}$ which cannot be $(w_f)$-stratified.

**2.6 Proposition.** *There exists a $C^p$ stratification of $f$.*

*Proof:*(c.f. [DM. Th. 4.8]) First note that if $f : \Gamma \longrightarrow \mathbf{R}^l$ is a $C^1$ definable map on a $C^1$-submanifold $\Gamma$ of $\mathbf{R}^n$, then the set

$$P = \{y \in \Gamma : \exists t > 0, \forall x \in \Gamma(\|x - y\| < t \Rightarrow \mathrm{rank} f(x) = \mathrm{rank} f(y))\}$$

is definable and $\dim(\Gamma \setminus P) < \dim \Gamma$.
Therefore, applying Proposition 1.4, we have a $C^p$ stratification of $f$. ◇

By the previous proposition and Proposition 1.4, Theorem 2.5 is implied by the following.

**2.7 Theorem.** *Suppose that $\mathcal{D}$ is polynomially bounded. Let $\Gamma, \Gamma'$ be definanable $C^p$ submanifolds of $\mathbf{R}^n$. Suppose $\Gamma \subset \overline{\Gamma'} \setminus \Gamma'$, and $f : \Gamma \cup \Gamma' \longrightarrow \mathbf{R}$ is continuous definable function such that $f$ has constant rank on $\Gamma$ and $\Gamma'$. Then*
*(i) $W_f = \{x \in \Gamma : (w_f) \text{ is satisfied at } x\}$ is definable, and*
*(ii) $\dim(\Gamma \setminus W_f) < \dim \Gamma$.*

*Proof:* The proof is much the same way as that for the condition $(a_f)$ in [L2].
(i) Since $x \mapsto d(f|_\Gamma)$ is a definable map (see [DM]), the kernel bundle of $f|_\Gamma$ is definable. Therefore,

$$\begin{aligned}W_f = \ &\{y_0 : y_0 \in \Gamma, \exists C > 0, \exists t > 0, \forall x \in \Gamma', \forall y \in \Gamma \\ &\|x - y_0\| < t, \|y - y_0\| < t \Rightarrow \delta(\ker d(f_\Gamma)(y), \ker d(f_{\Gamma'})(x) \leq C\|x - y\|\}\end{aligned}$$

is definable.
(ii) To prove the second assertion there are three cases to consider.
*Case 1:* $\mathrm{rank} f|_\Gamma = \mathrm{rank} f|_{\Gamma'} = 0$. In this case

$$W_f = \{y \in \Gamma : (\Gamma, \Gamma') \text{ satisfies Verdier condition (w) at } y\}.$$

The assertion follows from Theorem 1.3.
*Case 2:* $\mathrm{rank} f|_\Gamma = 0$ and $\mathrm{rank} f|_{\Gamma'} = 1$.
Suppose the contrary: $\dim(\Gamma \setminus W_f) < \dim \Gamma$. Since $(w_f)$ is $C^2$ invariant, by Cell Decomposition, we can assume that $\Gamma$ is an open subset of $\mathbf{R}^k \subset \mathbf{R}^k \times \mathbf{R}^{n-k}$, and $f|_{\Gamma'} > 0$, $f|_\Gamma \equiv 0$. So $T_{y,f} = \mathbf{R}^k$, for all $y \in \Gamma$. Let

$$A = \{(y, s, t) : (y, s) \in \Gamma \cup \Gamma', t > 0, f(y, s) = t\}.$$

Then $A$ is a definable set. By Definable Choice and the assumption, there exists an open subset $U$ of $\Gamma$, $\alpha > 0$, and a definable map $\theta : U \times [0, \alpha) \longrightarrow \mathbf{R}^{n-k}$, such that $\theta$ is $C^p$ on $U \times (0, \alpha)$, $\theta|_\Gamma \equiv 0$, and $f(y, \theta(y, t)) = t$, and, moreover, for all $y \in U$, we have

$$(*) \qquad \frac{\|D_1\theta(y,t)\|}{\|\theta(y,t)\|} \geq \frac{\delta(\mathbf{R}^k, T_{(y,\theta(y,t)),f})}{\|\theta(y,t)\|} \to +\infty, \quad \text{when } t \to 0^+.$$

On the other hand, by [M2. Prop. 5.2], there exists an open subset $B$ of $U$ and $r > 0$, such that

$$(**) \qquad \theta(y,t) = c(y)t^r + \varphi(y,t)t^{r_1}, y \in B, t > 0 \text{ sufficiently small,}$$

where $c$ is $C^p$ on $B$, $c \not\equiv 0$, $r_1 > r$, and $\varphi$ is $C^p$ with $\lim_{t \to 0^+} \varphi(y,t) = 0$, for all $y \in B$. Moreover, by Lemma 1.8, we can suppose that $D_1\varphi$ is bounded. Substituting $(**)$ to the left-hand side of $(*)$ we get a contradiction.

*Case 3:* $\text{rank} f|_\Gamma = \text{rank} f|_{\Gamma'} = 1$.
If $\dim(\Gamma \setminus W_f) = \dim \Gamma$, then the condition $(w_f)$ is false for $(\Gamma, \Gamma')$ over an open subset $B$ of $\Gamma$. It is easy to see that there is $c \in \mathbf{R}$ such that $(w_f)$ is false for the pair $(\Gamma \cap f^{-1}(c), \Gamma')$ over an open subset of $B \cap f^{-1}(c)$, and hence open in $\Gamma \cap f^{-1}(c)$. This contradicts Case 2. ◇

**2.8 Remark.** If the structure admits analytic cell decomposition, then the theorems hold true with "analytic" in place of "$C^p$". Our results can be translated to the setting of analytic-geometric categories in the sense of [DM].

DEPARTMENT OF MATHEMATICS
UNIVERSITY OF DALAT
DALAT - VIETNAM